\documentclass[12pt]{amsart}

\usepackage{amssymb}
\newtheorem{theorem}{Theorem}[section]
\newtheorem{lemma}[theorem]{Lemma}
\newtheorem{proposition}[theorem]{Proposition}
\newtheorem{remark}[theorem]{Remark}

\numberwithin{equation}{section}

\begin{document}

\baselineskip=15.5pt

\title[Torelli theorem for moduli of connections]{Torelli
theorem for the moduli spaces of connections on a Riemann surface}

\author[I. Biswas]{Indranil Biswas}

\address{School of Mathematics, Tata Institute of Fundamental
Research, Homi Bhabha Road, Bombay 400005, India}

\email{indranil@math.tifr.res.in}

\author[V. Mu{\~n}oz]{Vicente Mu{\~n}oz}

\address{Departamento de Matem\'aticas,
Consejo Superior de Investigaciones Cient{\'\i}ficas,
Serrano 113 bis, 28006 Madrid, Spain}

\email{vicente.munoz@imaff.cfmac.csic.es}

\subjclass[2000]{14D20, 14C34}

\keywords{Logarithmic connection, moduli space, Torelli theorem}

\date{}

\begin{abstract}
Let $(X\, , x_0)$ be any one--pointed compact
connected Riemann surface of genus $g$, with $g\geq 3$.
Fix two mutually coprime integers $r>1$ and $d$.
Let ${\mathcal M}_X$ denote the moduli space parametrizing all
logarithmic $\text{SL}(r,{\mathbb C})$--connections, singular
over $x_0$, on vector bundles over $X$ of degree $d$. We prove
that the isomorphism class of the variety ${\mathcal M}_X$
determines the Riemann surface $X$ uniquely up to an
isomorphism, although the biholomorphism class of ${\mathcal
M}_X$ is known to be independent of the complex structure of
$X$. The isomorphism class of the variety ${\mathcal M}_X$
is independent of the point $x_0\,\in\, X$. A similar
result is proved for the moduli space parametrizing
logarithmic $\text{GL}(r,{\mathbb C})$--connections,
singular over $x_0$, on vector bundles over $X$ of degree
$d$. The assumption $r>1$ is necessary for the moduli
space of logarithmic $\text{GL}(r,{\mathbb C})$--connections
to determine the isomorphism class of $X$ uniquely.
\end{abstract}

\maketitle

\section{Introduction}

Our aim here is to show that the algebraic
isomorphism class of a smooth moduli
space of $\text{SL}(r,{\mathbb C})$--connections on a
compact Riemann surface determines uniquely the isomorphism
class of the Riemann surface. This moduli space is
canonically biholomorphic to a certain
representation space of the fundamental group of the Riemann
surface; the biholomorphism is constructed by sending a connection
to its monodromy representation. Therefore, the biholomorphism
class of the moduli space is independent of the complex structure
of the Riemann surface. We note that any two smooth projective
varieties that are biholomorphic are actually
algebraically isomorphic.
However, this is not true for quasiprojective varieties.

Let $X$ be a compact connected Riemann surface of genus $g$,
with $g\, \geq\, 3$. Fix a base point $x_0\, \in\, X$. Take
any integer $r\, \geq\, 2$. Fix an integer $d$ which is
coprime to $r$. The de Rham differential $f\, \longmapsto\,
{\rm d}f$ defines a logarithmic connection
on the line bundle ${\mathcal O}_X(dx_0)$ which is singular
exactly over $x_0$.

Let ${\mathcal M}_X$ denote the moduli space
parametrizing logarithmic connections
$(E\, ,D)$ on $X$ singular exactly over $x_0$ and satisfying
the following three conditions:
\begin{itemize}
\item $E$ is a holomorphic vector bundle over $X$ of rank
$r$ with $\bigwedge^r E\, \cong\, {\mathcal O}_X(dx_0)$,

\item the logarithmic connection on $\bigwedge^r E$
induced by
the logarithmic connection $D$ on $E$
coincides with the connection on ${\mathcal O}_X(dx_0)$
given by the Rham differential, and

\item the residue $\text{Res}(D,x_0) \, =\,
-\frac{d}{r}\text{Id}_{E_{x_0}}$.
\end{itemize}
This moduli space ${\mathcal M}_X$ is an irreducible smooth
quasiprojective variety defined over $\mathbb C$ of dimension
$2(r^2-1)(g-1)$.
The isomorphism class of the variety ${\mathcal M}_X$ is
independent of the choice of the point $x_0\, \in\, X$;
see Remark \ref{rem1} for the details.

We prove the following Torelli type theorem:

\begin{theorem}\label{thm0}
The isomorphism class of the Riemann surface $X$ is uniquely
determined by the
isomorphism class of the variety ${\mathcal M}_X$. In other
words, if $(Y\, ,y_0)$ is another one--pointed
compact connected Riemann
surface of genus $g$ and ${\mathcal M}_Y$ the moduli space
of rank $r$ logarithmic connections $(E\, ,D)$ over $Y$
singular exactly over $y_0$ such that $\bigwedge^r E
\, \cong\,{\mathcal O}_X(dy_0)$, the
connection on $\bigwedge^r E$ induced by $D$ is the
de Rham connection on ${\mathcal O}_X(dy_0)$,
and ${\rm Res}(D,y_0) \, =\,-\frac{d}{r}{\rm Id}_{E_{y_0}}$,
then the two varieties ${\mathcal M}_X$ and
${\mathcal M}_Y$ are isomorphic if and only if the two Riemann
surfaces $X$ and $Y$ are isomorphic.
\end{theorem}

The proof of Theorem \ref{thm0} involves investigations of the
mixed Hodge structure of the moduli space ${\mathcal M}_X$.

In the special case where $r\, =\, 2\, < \, g-1$, a
similar result was proved in \cite{BN}.

Fix a point $x'\, \in\, X\setminus \{x_0\}\, =:\, X'$. Let
\[
\text{Hom}^0(\pi_1(X',x')\, , \text{SL}(r, {\mathbb C}))
\, \subset\, \text{Hom}(\pi_1(X',x')\, , \text{SL}(r, {\mathbb C}))
\]
be the space of homomorphisms from the fundamental
group $\pi_1(X',x')$ to $\text{SL}(r, {\mathbb C})$ satisfying the
condition that the image of the free homotopy class of oriented
loops in $X'$ around $x_0$ (with anticlockwise orientation) is
$\exp(2\pi\sqrt{-1}d/r)\cdot I_{r\times r}$. Any oriented loop in
$X'$ gives a conjugacy class in $\pi_1(X',x')$. The above condition
says that all elements in the conjugacy class are mapped
to $\exp(2\pi\sqrt{-1}d/r)\cdot I_{r\times r}$.

Let
\begin{equation}\label{rs}
{\mathcal R}_g\, :=\, \text{Hom}^0(\pi_1(X',x')\, , \text{SL}(r,
{\mathbb C}))/\text{SL}(r, {\mathbb C})
\end{equation}
be the quotient for the action of $\text{SL}(r, {\mathbb C})$
given by the
conjugation action of $\text{SL}(r, {\mathbb C})$ on itself.
The algebraic structure of $\text{SL}(r, {\mathbb C})$ induces
an algebraic structure on ${\mathcal R}_g$.
It is known that ${\mathcal R}_g$ is an irreducible smooth
quasiprojective variety defined over $\mathbb C$ of dimension
$2(r^2-1)(g-1)$. The variety ${\mathcal R}_g$ does not depend on the
choice of the point $x'$ as the fundamental groups of a space
for two different
base points are identified up to an inner automorphism. The
isomorphism class of the variety
${\mathcal R}_g$ is independent of the
complex structure of $X$. It depends only on the genus of $X$
for fixed $r$ and $d$.

Given any logarithmic connection $(E\, ,D)\, \in\, {\mathcal M}_X$,
consider the corresponding monodromy
representation of $\pi_1(X',x')$ on the fiber
$E_{x'}$. From the conditions on $D$ it follows that this
monodromy representation lies in
$\text{Hom}^0(\pi_1(X',x')\, , \text{SL}(r,{\mathbb C}))$.
Consequently, we get a map
\[
\sigma\, :\, {\mathcal M}_X\, \longrightarrow\, {\mathcal R}_g
\]
that sends any connection to its monodromy representation.
It is known that $\sigma$ is actually a biholomorphism
\cite[p. 26, Theorem 7.1]{Si3}.
This immediately implies that the complex manifold
${\mathcal M}_X$ is biholomorphic to the
complex manifold ${\mathcal M}_Y$ defined in Theorem
\ref{thm0}.

Take $r$ and $d$ as above.
Let $\widehat{\mathcal M}_X$ denote the moduli space
parametrizing logarithmic connections $(E\, , D)$
on $X$, singular exactly over $x_0$, such that
\begin{itemize}
\item $\text{rank}(E) \, =\, r$, and

\item $\text{Res}(D,x_0) \, =\,
-\frac{d}{r}\text{Id}_{E_{x_0}}$.
\end{itemize}
The second condition implies that $\text{degree}(E)\, =\, d$.
It is known that $\widehat{\mathcal M}_X$ is a smooth
irreducible quasiprojective variety defined over $\mathbb C$
of dimension $2(r^2(g-1)+1)$.
The isomorphism class of the variety $\widehat{\mathcal M}_X$
is independent of the choice of the point $x_0$.

We prove the following analog of Theorem \ref{thm0}:

\begin{theorem}\label{thm00}
The isomorphism class of the Riemann surface $X$ is uniquely
determined by the
isomorphism class of the variety $\widehat{\mathcal M}_X$.
\end{theorem}

Since $\widehat{\mathcal M}_X$ is biholomorphic to the
representation space of the fundamental
group obtained by replacing
$\text{SL}(n,{\mathbb C})$ with $\text{GL}(n,{\mathbb C})$
in Eqn. \eqref{rs}, the biholomorphism class of
$\widehat{\mathcal M}_X$ is also independent of the complex
structure of $X$.

It should be mentioned that Theorem \ref{thm00} fails for
$r\, =\, 1$; see Remark \ref{ra1} for the details.

We do not know if Theorem \ref{thm0} and Theorem \ref{thm00}
remain valid when $g\, =\,2$. In our proof of
Theorem \ref{thm0} and Theorem \ref{thm00}, we needed the
codimension in Lemma \ref{lem.1} to be at least two.
We should mention that this is
the only place where our assumption $g\, > \, 2$ is used.

\section{Biholomorphism class of the moduli space}\label{sec1}

Let $X$ be a compact connected Riemann surface. The genus
of $X$, which will be denoted by $g$, is assumed to be
at least three. The holomorphic cotangent bundle of $X$ will
be denoted by $K_X$. The holomorphic tangent bundle of $X$ will
be denoted by $TX$. Fix a base point $x_0\, \in\, X$.

Take a holomorphic vector bundle $E$ over $X$.
A \textit{logarithmic connection on $E$ singular over $x_0$}
is a holomorphic differential operator
\begin{equation}\label{deD}
D\, :\, E\, \longrightarrow \, E\otimes K_X\otimes
{\mathcal O}_X(x_0)
\end{equation}
satisfying the Leibniz identity which says that
\[
D(fs) \, =\, f\cdot D(s) + s\otimes {\rm d}f\, ,
\]
where $f$ is any locally defined holomorphic function
on $X$, and
$s$ is any locally defined holomorphic section of $E$. The
above Leibniz identity implies that the order of $D$ is
exactly one. More precisely, a logarithmic connection on $E$
singular over $x_0$ is a holomorphic differential operator
of order one
\[
D\, \in\, H^0(X,\, \text{Diff}^1_X(E,E\otimes K_X\otimes
{\mathcal O}_X(x_0)))
\]
whose symbol, which is a holomorphic section of the vector
bundle
\[
\text{Hom}(E,E\otimes K_X\otimes
{\mathcal O}_X(x_0))\otimes TX \, =\,
\text{End}(E)\otimes_{{\mathcal O}_X}
{\mathcal O}_X(x_0)\, ,
\]
coincides with the section of $\text{End}(E)
\bigotimes_{{\mathcal O}_X}{\mathcal O}_X(x_0)$
given by the identity automorphism of the vector bundle $E$.

The fiber over $x_0$ of the line bundle
$K_X\bigotimes_{{\mathcal O}_X}{\mathcal O}_X(x_0)$
is canonically identified
with $\mathbb C$. This identification is obtained by sending
any holomorphic section of
$K_X\bigotimes_{{\mathcal O}_X}{\mathcal O}_X(x_0)$
defined around $x_0$,
which is a meromorphic form with a pole at $x_0$
of order at most one, to its residue at $x_0$.

Let $D$ be a logarithmic connection on $E$
singular over $x_0$. Consider the composition
\[
E \, \stackrel{D}{\longrightarrow}\, E\otimes
K_X\otimes {\mathcal O}_X(x_0) \, \longrightarrow\,
(E\otimes K_X\otimes {\mathcal O}_X(x_0))_{x_0}
\, =\, E_{x_0}\, ,
\]
where the last equality is obtained using the above mentioned
identification of the line
$(K_X\bigotimes {\mathcal O}_X(x_0))_{x_0}$
with $\mathbb C$, and the second homomorphism is the evaluation
at $x_0$ of sections. This composite homomorphism is
${\mathcal O}_X$--linear. Hence it gives a linear endomorphism
of the fiber $E_{x_0}$. This endomorphism is called the
\textit{residue} of $D$ at $x_0$, and it is denoted by
$\text{Res}(D,x_0)$. (See \cite[p. 53]{De1}.) We have
\begin{equation}\label{de-re}
\text{degree}(E) \, =\, -\text{trace}(\text{Res}(D,x_0))
\end{equation}
(see the last sentence of Corollary B.3 in \cite[p. 186]{EV}).

Fix an
integer $r \, \geq\, 2$, and also fix an integer $d$ which is
coprime to $r$. Let $D_0$ denote the logarithmic connection
on the line bundle ${\mathcal O}_X(dx_0)$, singular over $x_0$,
defined by the de Rham differential that sends any locally
defined meromorphic function $f$ on $X$ to the form ${\rm d}f$;
here we identify the sheaf of sections of the line bundle
${\mathcal O}_X(dx_0)$ with the sheaf of meromorphic
functions on $X$ with a pole at $x_0$ of order at most $d$. It
is straight--forward to check that the residue of $D_0$ at $x_0$
is $-d$. Note that this also follows from Eqn. \eqref{de-re}.

Let ${\mathcal M}_X$ denote the moduli space parametrizing
all pairs $(E\, , D)$ of the following type:
\begin{itemize}
\item $E$ is a holomorphic vector bundle
of rank $r$ over $X$ with
$\bigwedge^r E \, \cong\, {\mathcal O}_X(dx_0)$,

\item $D$ is a logarithmic connection on $E$ singular
over $x_0$ with
\[
\text{Res}(D,x_0) \, =\,-\frac{d}{r}\text{Id}_{E_{x_0}}
\]

\item the connection on $\bigwedge^r E \, \cong\,
{\mathcal O}_X(dx_0)$ induced by $D$ coincides, by
some (hence any) isomorphism between $\bigwedge^r E$
and ${\mathcal O}_X(dx_0)$, with the de Rham connection
$D_0$ on ${\mathcal O}_X(dx_0)$ defined above.
\end{itemize}
See \cite{Si2}, \cite{Si3} and \cite{Ni2} for the construction
of this moduli space ${\mathcal M}_X$. The scheme
${\mathcal M}_X$ is a reduced and irreducible
quasiprojective variety defined over $\mathbb C$, and its (complex)
dimension is $2(r^2-1)(g-1)$, where $g$ is the genus of $X$.

Take any logarithmic connection $(E\, , D)\, \in\,
{\mathcal M}_X$. Since $d$ and $r$ are mutually
coprime, the connection $D$ is irreducible in the sense that
no nonzero proper holomorphic subbundle of $E$ is preserved
by $D$ \cite[p. 787, Lemma 2.3]{BR}. Using this it
follows that the variety ${\mathcal M}_X$ is smooth. See
\cite{Ha} for some topological properties of ${\mathcal M}_X$.

Let $X' \, :=\, X \setminus \{x_0\}$
be the complement. Fix a point $x'\, \in\, X'$.
The point $x_0$ gives a conjugacy class in the fundamental group
$\pi_1(X',x')$ as follows: Let
\[
f\, :\, {\mathbb D} \, \longrightarrow\, X'
\]
be an orientation preserving embedding of the closed unit
disk ${\mathbb D}\, :=\, \{z\,\in\,
{\mathbb C}\, \vert\, ||z||^2 \,\leq\, 1\}$
into the Riemann surface $X$ such that $f(0) \, =\, x_0$. The
free homotopy class of the map $\partial {\mathbb D}
\, =\, S^1 \, \longrightarrow\, X'$ obtained by
restricting $f$ to the boundary of ${\mathbb D}$ is independent
of the choice of $f$. The orientation of $\partial {\mathbb D}$
coincides with the anti--clockwise rotation around $x_0$. Any
free homotopy class of oriented loops in $X'$
gives a conjugacy class in $\pi_1(X',x')$. Let $\gamma$ denote the
orbit in $\pi_1(X',x')$, for the conjugation action of
$\pi_1(X',x')$ on itself, defined by the above free
homotopy class of oriented loops associated to $x_0$.

Let
\begin{equation}\label{s.r.}
\text{Hom}^0(\pi_1(X',x')\, , \text{SL}(r, {\mathbb C}))
\, \subset\, \text{Hom}(\pi_1(X',x')\, , \text{SL}(r, {\mathbb C}))
\end{equation}
be the space of homomorphisms from $\pi_1(X',x')$ to
$\text{SL}(r, {\mathbb C})$ satisfying the condition that the
image of $\gamma$ (defined above) is
$\exp(2\pi\sqrt{-1}d/r)\cdot I_{r\times r}$. It should be clarified that
since $\exp(2\pi\sqrt{-1}d/r)\cdot I_{r\times r}$ is in the center of
$\text{SL}(r, {\mathbb C})$, a homomorphism sends the orbit
$\gamma$ in $\pi_1(X',x')$ (for the adjoint
action of $\pi_1(X',x')$ on itself)
to $\exp(2\pi\sqrt{-1}d/r)\cdot I_{r\times r}$ if and only if there
is an element in the orbit which is mapped to
$\exp(2\pi\sqrt{-1}d/r)\cdot I_{r\times r}$.

Take any homomorphism $\rho \, \in\, \text{Hom}^0(\pi_1(X',x')\, ,
\text{SL}(r, {\mathbb C}))$. Let
$(V\, , \nabla)$ be the flat vector bundle or rank $r$
over $X'$ given by $\rho$. The monodromy of $\nabla$ along the
oriented loop
$\gamma$ is $\exp(2\pi\sqrt{-1}d/r)\cdot I_{r\times r}$.
Using the logarithm
$2\pi\sqrt{-1}d/r\cdot I_{r\times r}$ of the
monodromy, the vector bundle $V$ over $X'$ extends to a holomorphic
vector bundle $\overline{V}$ over $X$, and furthermore, the connection
$\nabla$ on $V$ extends to a logarithmic connection
$\overline{\nabla}$ on the vector
bundle $\overline{V}$ over $X$ such that
$(\overline{V}\, ,\overline{\nabla})\, \in\, {\mathcal M}_X$
(see \cite[p. 159, Theorem 4.4]{Ma}).

The easiest way to construct the
above mentioned extension $(\overline{V}\, ,\overline{\nabla})$
is to fix a ramified Galois covering
\begin{equation}\label{g.co.}
\phi\, :\, Y\, \longrightarrow\, X
\end{equation}
of degree $r$ which is totally ramified
over $x_0$ (it may have other ramification
points). Let $y_0\,\in\, Y$ be the unique point
such that
\begin{equation}\label{pt.}
\phi(y_0)\, =\, x_0\, .
\end{equation}
The connection $\phi^*\nabla$ on the vector bundle
$\phi^*V$ over $\phi^{-1}(X')\,=\, Y\setminus\{y_0\}
\, \subset\, Y$
has trivial monodromy around the point $y_0$.
Indeed, this follows from the given condition that the
map $\phi$ is totally ramified over $x_0$, combined with
the fact that the order of
the matrix $\exp(2\pi\sqrt{-1}d/r)\cdot I_{r\times r}$ is a
divisor of the degree of $\phi$.
Consequently, the flat vector bundle
$(\phi^*V\, , \phi^*\nabla)$
has a canonical extension to $Y$ as a flat
vector bundle (take the direct image of the local
system by the inclusion map $\phi^{-1}(X')\,\hookrightarrow
\, Y$). Let $(W\, , \nabla')$ denote this
flat vector bundle over $Y$ obtained by extending
$(\phi^*V\, , \phi^*\nabla)$.
Let $\nabla_0$ denote the de Rham logarithmic connection
on the line bundle ${\mathcal O}_Y(dy_0)$ that sends
any locally defined meromorphic function $f$ to ${\rm d}f$.
The vector bundle $W\bigotimes {\mathcal O}_Y(dy_0)$
is equipped with the logarithmic connection
\[
\nabla''\, :=\, \nabla'\otimes \text{Id}_{{\mathcal O}_Y(dy_0)}
+ \text{Id}_W \otimes\nabla_0\, .
\]
The direct image $\phi_* (W\bigotimes {\mathcal O}_Y(dy_0))$
over $X$ is equipped with an action of the Galois group
$\text{Gal}(\phi)$ for $\phi$, with $\text{Gal}(\phi)$
acting trivially on $X$.
The above mentioned holomorphic vector bundle
$\overline{V}$ over $X$ is identified with the invariant part
$(\phi_* (W\bigotimes {\mathcal O}_Y(dy_0)))^{\text{Gal}(\phi)}$
for this action of $\text{Gal}(\phi)$
on $\phi_* (W\bigotimes {\mathcal O}_Y(dy_0))$.
The logarithmic connection $\overline{\nabla}$ on
$\overline{V}$ coincides with the one induced by the
logarithmic connection $\nabla''$ on $W\bigotimes {\mathcal O}_Y(dy_0)$.

Since $\text{SL}(r, {\mathbb C})$ is an algebraic group
defined over $\mathbb C$, and $\pi_1(X',x')$ is a
finitely presented group, the representation
space $\text{Hom}(\pi_1(X',x')\, , \text{SL}(r, {\mathbb C}))$
is a complex algebraic
variety in a natural way. The closed subvariety
$\text{Hom}^0(\pi_1(X',x')\, , \text{SL}(r, {\mathbb C}))$
defined in Eqn. \eqref{s.r.} is actually smooth. The
smoothness of $\text{Hom}^0(\pi_1(X',x')\, , \text{SL}(r, {\mathbb
C}))$ follows from the fact that any representation in
$\text{Hom}^0(\pi_1(X',x')\, , \text{SL}(r, {\mathbb C}))$
is irreducible; the irreducibility of such a
representation follows from \cite[p. 787, Lemma 2.3]{BR}.
The conjugation action of $\text{SL}(r, {\mathbb C})$
on itself induces an action of $\text{SL}(r, {\mathbb C})$ on
$\text{Hom}^0(\pi_1(X',x')\, , \text{SL}(r, {\mathbb C}))$. The
action of any $T\, \in\, \text{SL}(r, {\mathbb C})$ on
$\text{Hom}^0(\pi_1(X',x')\, , \text{SL}(r, {\mathbb C}))$ sends
any homomorphism
$\rho$ to the homomorphism $\pi_1(X',x')\, \longrightarrow\,
\text{SL}(r, {\mathbb C})$ defined by $\beta\, \longmapsto\,
T^{-1}\rho(\beta)T$. Let
\begin{equation}\label{de.R}
{\mathcal R}_g\, :=\, \text{Hom}^0(\pi_1(X',x')\, , \text{SL}(r,
{\mathbb C}))/\text{SL}(r, {\mathbb C})
\end{equation}
be the quotient space for this action.

The algebraic structure of $\text{Hom}^0(\pi_1(X',x')\, , \text{SL}(r,
{\mathbb C}))$ induces an algebraic structure on the quotient
${\mathcal R}_g$. The scheme ${\mathcal R}_g$ is an irreducible smooth
quasiprojective variety of dimension $2(r^2-1)(g-1)$ defined
over $\mathbb C$. We recall that ${\mathcal R}_g$ and ${\mathcal M}_X$
are known as the \textit{Betti moduli space} and the
\textit{de Rham moduli space} respectively (see \cite{Si1},
\cite{Si2}, \cite{Ha}).

If we replace the base point $x'\, \in\, X'$ by another point
$x''\, \in\, X'$, then we can construct an isomorphism of
$\pi_1(X',x')$ with $\pi_1(X',x'')$ by fixing a path
connecting $x'$ to $x''$. Consequently, $\pi_1(X',x')$
and $\pi_1(X',x'')$ are identified up to an inner automorphism.
This in turn implies that the variety ${\mathcal R}_g$ is canonically
identified with the variety obtained by replacing $x'$ with $x''$
in the construction of ${\mathcal R}_g$. In other words,
${\mathcal R}_g$ does not depend on the choice of the point $x'$.

Given any two one--pointed compact connected oriented
$C^\infty$ surfaces
of genus $g$, there is an orientation preserving diffeomorphism
between them
that takes the marked point in one surface to the
marked point in the other surface. Therefore, the isomorphism
class of the variety ${\mathcal R}_g$ depends only on
the integers
$g$, $r$ and $d$. In particular, the isomorphism class
of this variety is independent of the complex structure of the
topological surface $X$.

Given any $(E\, ,D)\,\in\, {\mathcal M}_X$, the monodromy
representation of $D$
\[
\pi_1(X',x')\, \longrightarrow\, \text{Aut}(E_{x'})
\]
gives an element of ${\mathcal R}_g$. The inverse map is
obtained from the earlier mentioned construction that associates
an element of ${\mathcal M}_X$ to each element of
$\text{Hom}^0(\pi_1(X',x')\, , \text{SL}(r, {\mathbb C}))$.

These two maps are inverses of each other \cite[p. 26, Theorem
7.1]{Si3}. In particular, the following proposition holds.

\begin{proposition}\label{prop.1}
The moduli space ${\mathcal M}_X$ is biholomorphic to
${\mathcal R}_g$. Therefore, if $(Y\, ,y_0)$ is another
one--pointed compact connected Riemann surface of genus $g$
and ${\mathcal M}_Y$ the corresponding moduli space of
logarithmic connections, then the two complex manifolds
${\mathcal M}_X$ and ${\mathcal M}_Y$ are biholomorphic.
\end{proposition}

\begin{remark}\label{rem1}
{\rm The isomorphism class of the variety ${\mathcal M}_X$
is independent of the choice of the point $x_0\, \in\, X$. To
prove this, take another point $x_1\, \in\, X$. Let
${\mathcal M}^{x_1}_X$ be the moduli space of
logarithmic connections
on $X$, singular exactly over $x_1$, obtained by replacing
$x_0$ with $x_1$ in the construction of ${\mathcal M}_X$.

Fix a holomorphic line bundle $L$ over $X$ such that
${\mathcal O}_X(x_1-x_0)$ is isomorphic to $L^{\otimes d}$.
Let $D'_0$ denote the de Rham logarithmic
connection on the line bundle ${\mathcal O}_X(x_1-x_0)$
that sends any locally defined meromorphic function $f$ to
${\rm d}f$.
We note that there is a unique logarithmic connection $D_0$
on $L$, singular over both $x_0$ and $x_1$, such that the
logarithmic connection on the tensor product
$L^{\otimes d}\, \cong\,{\mathcal O}_X(x_1-x_0)$ induced
by $D_0$ coincides with the de Rham connection $D'_0$ on
${\mathcal O}_X(x_1-x_0)$.

For any $(E\, ,D)\, \in\, {\mathcal M}_X$, it is easy to
see that
$$
(E\otimes L\, , D\otimes \text{Id}_{L}+\text{Id}_E
\otimes D_0)\, \in\, {\mathcal M}^{x_1}_X\, .
$$
The map
$$
{\mathcal M}_X\, \longrightarrow\,
{\mathcal M}^{x_1}_X
$$
defined by
$$
(E\, ,D)\, \longmapsto\,
(E\otimes L\, , D\otimes \text{Id}_{L}+\text{Id}_E
\otimes D_0)
$$
is an algebraic isomorphism of varieties.}
\end{remark}

In the next section we will investigate the algebraic
structure of ${\mathcal M}_X$.

\section{The second intermediate Jacobian of the moduli
space}\label{sec3}

Let
\begin{equation}\label{cU}
{\mathcal U}\, \subset\, {\mathcal M}_X
\end{equation}
be the Zariski open
subset parametrizing all $(E\, ,D)$ such that the underlying
vector bundle $E$ is stable. The openness of this subset
follows from
\cite[p. 182, Proposition 10]{Sh}. Let ${\mathcal N}_X$ denote the
moduli space parametrizing all stable vector bundles $E$ over $X$
with $\text{rank}(E)\,=\, r$ and $\bigwedge^r E\,\cong\,{\mathcal
O}_X(dx_0)$. The moduli space ${\mathcal N}_X$ is an irreducible
smooth projective variety of dimension $(r^2-1)(g-1)$
defined over $\mathbb C$. Let
\begin{equation}\label{Ph}
\Phi\, :\, {\mathcal U}\, \longrightarrow\, {\mathcal N}_X
\end{equation}
denote the forgetful map that sends any $(E\, ,D)$ to $E$.

It is know that any $E\, \in\, {\mathcal N}_X$ admits a logarithmic
connection $D$ such that $(E\, ,D)\, \in\, {\mathcal M}_X$
\cite[p. 787, Lemma 2.3]{BR}. Therefore, the projection $\Phi$
in Eqn. \eqref{Ph} is surjective. Furthermore, $\Phi$
makes ${\mathcal U}$ an affine bundle over
${\mathcal N}_X$. More precisely, ${\mathcal U}$ is a
torsor over ${\mathcal N}_X$
for the holomorphic cotangent bundle $T^*{\mathcal N}_X$.
This means that the fibers of the vector
bundle $T^*{\mathcal N}_X$ act freely
transitively on the fibers of $\Phi$ \cite[p. 786]{BR}.
Since ${\mathcal M}_X$ is irreducible, and ${\mathcal U}$
is nonempty, the subset ${\mathcal U}\,\subset\,
{\mathcal M}_X$ is Zariski dense.

\begin{lemma}\label{lem.1}
Let ${\mathcal Z}\, :=\, {\mathcal M}_X\setminus
{\mathcal U}$ be the complement of the Zariski open
dense subset. The codimension of the Zariski
closed subset ${\mathcal Z}$ in ${\mathcal M}_X$ is at least
$(r-1)(g-2) +1$.
\end{lemma}

\begin{proof}
Let $E$ be a holomorphic vector bundle over $X$
with $\text{rank}(E)\,=\, r$ and
$\bigwedge^r E\, \cong\, {\mathcal O}_X(dx_0)$. Assume
that $E$ admits a logarithmic connection $D$ singular
over $x_0$ such that $(E\, , D)\, \in\, {\mathcal M}_X$.

Fix a logarithmic connection $D$ on $E$ singular
over $x_0$ such that $(E\, , D)\, \in\, {\mathcal M}_X$.
Let $T\, \in\, H^0(X,\, \text{End}(E))$ be a holomorphic
endomorphism of $E$ satisfying the identity
\begin{equation}\label{Id.}
D\circ T\, =\, (T\otimes \text{Id}_{K_X\otimes
{\mathcal O}_X(x_0)})\circ D\, .
\end{equation}
We will show that there is a complex number $\lambda$
such that $T\, =\, \lambda\cdot \text{Id}_{E}$, or in other
words, any endomorphism of the pair $(E\, ,D)$ is a scalar
multiplication.

To prove this, first note that as $X$ is a compact
connected Riemann surface, it does not admit any nonconstant
holomorphic functions. In particular, $\text{trace}(T^j)$
is a constant function on $X$ for all $j\, \geq\, 1$. This
implies that the eigenvalues, as well as their multiplicities,
of $T(x)\, \in\, \text{End}(E_x)$ are independent of
$x\, \in\, X$. Consequently, for any eigenvalue
$\alpha$ of $T(x)$, where $x$ is some point of $X$,
the generalized eigenspaces of the endomorphisms
$T(y)$, $y\, \in\, X$,
for the eigenvalue $\alpha$ patch together to
define a holomorphic subbundle of $E$
of positive rank. Let $E^\alpha$ denote this holomorphic
subbundle of $E$ defined by the generalized eigenspaces
for the eigenvalue $\alpha$.

{}From the identity in Eqn. \eqref{Id.} it follows immediately
that $E^\alpha$ is preserved by the connection $D$. Since
$D$ is irreducible \cite[p. 787, Lemma 2.3]{BR}, this implies
that $E^\alpha\, =\, E$. In other words, for any $y\,\in\, X$,
the endomorphism $T(y) \, \in\, \text{End}(E_y)$
has exactly one eigenvalue, namely $\alpha$.

Consider the endomorphism
\[
T_0\, =\, T - \alpha\cdot \text{Id}_{E}\, \in\,
H^0(X,\, \text{End}(E))\, .
\]
Since $\alpha$ is the only eigenvalue of $T(y)$, $y\, \in
\, X$, it follows that $T_0$ is a nilpotent endomorphism
of $E$. From Eqn. \eqref{Id.} it follows immediately that
\begin{equation}\label{Id2.}
D\circ T_0\, =\, (T_0\otimes \text{Id}_{K_X\otimes
{\mathcal O}_X(x_0)})\circ D\, .
\end{equation}

Let $E_0\, \subset\, E$ be the coherent subsheaf defined by
the kernel of the above endomorphism $T_0$. Note that as
$T_0$ is nilpotent, the subsheaf $E_0$ is
nonzero. From Eqn. \eqref{Id2.}
it follows that $E_0$ is preserved by the connection $D$.
Since $D$ is irreducible, this implies that $E_0\, =\, E$.
Consequently, $T_0\, =\, 0$. In other words,
we have $T \,=\, \alpha\cdot \text{Id}_{E}$.

Let ${\mathcal D}$ denote the space of all logarithmic
connections $D$ on $E$, singular exactly over the point
$x_0$, such that $(E\, , D)\, \in\, {\mathcal M}_X$.
We note that ${\mathcal D}$ is an affine space for the vector
space $H^0(X,\, \text{ad}(E)\bigotimes K_X)$, where
$\text{ad}(E)\, \subset\, \text{End}(E)$ is the subbundle
of rank $r^2-1$ defined by the sheaf
of endomorphisms of $E$ of trace zero
\cite[p. 786]{BR}. We have shown above that for the
natural action of the global automorphism group
$\text{Aut}(E)$ on ${\mathcal D}$, the isotropy at any point
is the subgroup defined by all automorphisms of the form
$\lambda\cdot \text{Id}_{E}$ with $\lambda\, \in\,{\mathbb C}^*$.

Therefore, the space of all isomorphism classes of logarithmic
connections $D$, on the given vector bundle $E$, such that
$(E\, ,D)\, \in\, {\mathcal M}_X$ is of dimension
\[
h^0(X,\, \text{ad}(E)\otimes K_X) -
(\dim \text{Aut}(E) - 1) \, =\,
h^1(X,\, \text{ad}(E)) - h^0(X,\, \text{ad}(E))
\, =\, (r^2-1)(g-1)\, ,
\]
where the last equality is the Riemann--Roch formula, and the
first equality follows from the Serre duality.

Therefore, if a holomorphic vector bundle $E'$ over $X$ admits
a logarithmic connection $D'$ singular
over $x_0$ such that $(E'\, ,D')\, \in\,
{\mathcal M}_X$, then the dimension of the space of all
isomorphism classes of such logarithmic connections on $E'$
is actually independent of the choice of $E'$.

We will now estimate the dimension of the isomorphism classes
of non--semistable vector bundles over $X$ that arise in the
family parametrized by ${\mathcal M}_X$.

Take any $(E\, ,D)\, \in\, {\mathcal M}_X$ such that the
underlying vector bundle $E$ is not stable. Since $r$ and $d$ are
mutually coprime, this implies that $E$ is not semistable. Let
\begin{equation}\label{filt.}
0\, =\, E_0 \, \subset\, E_1 \, \subset\, E_2
 \, \subset\, \cdots\, \subset\, E_{\ell-1}
 \, \subset\, E_{\ell} \, =\, E
\end{equation}
be the Harder--Narasimhan filtration of $E$.
We recall that the collection of pairs of integers
$\{(\text{rank}(E_i)\, ,
\text{degree}(E_i))\}_{i=1}^\ell$ is called the
\textit{Harder--Narasimhan polygon} of $E$ (see
\cite[p. 173]{Sh}).

The space of all isomorphism classes of holomorphic vector
bundles over $X$, whose Harder--Narasimhan polygon coincides
with that of the given vector bundle $E$, is of dimension at
most $r^2(g-1)-(r-1)(g-2)$. This follows from
\cite[p. 247--248]{Bh}; we will give below the details of the
argument.

Since
\[
\frac{\text{degree}(E_i/E_{i-1})}{\text{rank}(E_i/E_{i-1})}
\, >\, \frac{\text{degree}(E_{i+1}/E_{i})}{\text{rank}(E_{i+1}/E_i)}
\]
for all $i\, \in\, [1\, , \ell-1]$ in Eqn. \eqref{filt.},
the integer $c$
in the middle of page 247 of \cite{Bh} can be taken to be zero.
This follows immediately from the inductive construction of $c$
in \cite[p. 248]{Bh}
(see the comment in \cite[p. 247]{Bh} just after $c$ is
introduced). In other words, $a_r$ defined in
\cite[p. 247]{Bh} satisfies the inequality
\[
a_r \, \leq\, r^2(g-1)-(r-1)(g-2)
\]
(see \cite[p. 247]{Bh}); note that the term
$\sum_{1\leq i< j\leq r} n_in_j$ in
\cite[p. 247]{Bh} is at least $(\sum_i n_i)-1$
(the index $\ell$ here is $r$ in \cite{Bh}, and
$\sum_i n_i$ in \cite{Bh} is $r$ here).
The dimension of the space of all
isomorphism classes of vector bundles over $X$ whose
Harder--Narasimhan polygon coincides with that of $E$ is at
most $a_r$ \cite[p. 248]{Bh}.

Therefore, the space of all isomorphism classes of
vector bundles $E'$ over $X$ whose
Harder--Narasimhan polygon coincides with that of $E$,
and $\bigwedge^r E' \, \cong\, {\mathcal O}_X(dx_0)$,
is of dimension at most $r^2(g-1)-(r-1)(g-2)-g$.

Since the dimension of the space of all isomorphism
classes of logarithmic connections, lying in ${\mathcal M}_X$,
on any given vector bundle $E'$ is $(r^2-1)(g-1)$
(assuming that it admits a connection), the subvariety
of ${\mathcal M}_X$ parametrizing all pairs of the form
$(E'\, , D')\, \in\, {\mathcal M}_X$ such that the
Harder--Narasimhan polygon
of $E'$ coincides with that of $E$ is at most
$r^2(g-1)-(r-1)(g-2)-g+ (r^2-1)(g-1)$; that this subset
of ${\mathcal M}_X$ is algebraic follows from \cite[p. 182,
Proposition 10]{Sh}. We also note
that there are only finitely many
Harder--Narasimhan polygons that occur for the
vector bundles over $X$ in the family
parametrized by ${\mathcal M}_X$
\cite[p. 183, Proposition 11]{Sh}.

Since $\dim {\mathcal M}_X\, =\, 2(r^2-1)(g-1)$, we have
\[
\dim {\mathcal M}_X -(r^2(g-1)-(r-1)(g-2)-g+(r^2-1)(g-1))
\, =\, (r-1)(g-2)+1\, .
\]
This completes the proof of the lemma.
\end{proof}

For any $i\, \geq\, 0$, the $i$-th cohomology of a complex
variety with coefficients in $\mathbb Z$ is
equipped with a mixed Hodge structure \cite{De2}, \cite{De3}.

\begin{proposition}\label{prop2}
The mixed Hodge structure $H^3({\mathcal M}_X,\,{\mathbb Z})$
is pure of weight three. More precisely, the mixed Hodge
structure $H^3({\mathcal M}_X,\, {\mathbb Z})$ is isomorphic
to the Hodge structure $H^3({\mathcal N}_X,\, {\mathbb Z})$,
where ${\mathcal N}_X$ is the moduli space of stable vector
bundles introduced at the beginning of this section.
\end{proposition}

\begin{proof}
Consider the diagram of morphisms
\begin{equation}\label{diag.}
{\mathcal N}_X \, \stackrel{\Phi}{\longleftarrow}\,
{\mathcal U}\, \stackrel{\iota}{\longrightarrow}\,
{\mathcal M}_X\, ,
\end{equation}
where $\Phi$ is the projection in Eqn. \eqref{Ph},
and $\iota$ is the inclusion map.

Since ${\mathcal U}\,\stackrel{\Phi}{\longrightarrow}\,
{\mathcal N}_X$ is a fiber bundle with contractible fibers,
the induced homomorphism
\begin{equation}\label{is.}
\Phi^* \, :\, H^i({\mathcal N}_X,\, {\mathbb Z})\,
\longrightarrow\, H^i({\mathcal U},\, {\mathbb Z})
\end{equation}
is an isomorphism for all $i\, \geq\, 0$.
Therefore, $\Phi$ induces an isomorphism of the
two mixed Hodge structures $H^i({\mathcal N}_X,\, {\mathbb Z})$
and $H^i({\mathcal U},\, {\mathbb Z})$. Since ${\mathcal N}_X$
is a smooth projective variety, we know
that $H^i({\mathcal N}_X,\, {\mathbb Z})$ is a
pure Hodge structure of weight $i$. Thus $H^i({\mathcal U},\,
{\mathbb Z})$ is a pure Hodge structure of weight $i$. In
particular, $H^3({\mathcal U},\, {\mathbb Z})$ is a pure Hodge
structure of weight three.

Let
\begin{equation}\label{ex.m.}
H^3({\mathcal M}_X\, , {\mathcal U}\, ,{\mathbb Z}) \,
\longrightarrow\, H^3({\mathcal M}_X\, ,{\mathbb Z})\,
\stackrel{\psi}{\longrightarrow}\, H^3({\mathcal U}\, ,{\mathbb Z})
\,\longrightarrow\,
H^4({\mathcal M}_X\, , {\mathcal U}\, ,{\mathbb Z})
\end{equation}
be the long exact sequence of relative cohomologies. We note
that the above homomorphism $\psi$ is a morphism
of mixed Hodge structures. In fact, Eqn. \eqref{ex.m.}
itself is an exact sequence of mixed Hodge structures
\cite[p. 43, Proposition (8.3.9)]{De3}, but we will not need
this here. Let
\[
{\mathcal Z}\, :=\, {\mathcal M}_X \setminus {\mathcal U}
\]
be the complement, which is a closed subscheme.
Since $g\, \geq\, 3$ and $r\, \geq\, 2$, from Lemma \ref{lem.1}
we know that the (complex) codimension of ${\mathcal Z}$ in
${\mathcal M}_X$ is at least two. Therefore, the relative cohomology
has the following properties:
\begin{equation}\label{co.va.}
H^i({\mathcal M}_X\, ,{\mathcal U}\, ,{\mathbb Z})\, =\, 0
\end{equation}
for all $i\, <\, 4$, and $H^4({\mathcal M}_X\, ,
{\mathcal U}\, ,{\mathbb Z})$ is torsionfree. Consequently,
Eqn. \eqref{ex.m.} gives the exact sequence
\begin{equation}\label{ex.m.2}
0 \, \longrightarrow\, H^3({\mathcal M}_X\, ,{\mathbb Z})\,
\stackrel{\psi}{\longrightarrow}\, H^3({\mathcal U}\, ,{\mathbb Z})
\, \longrightarrow\, H^4({\mathcal M}_X\, ,
{\mathcal U}\, ,{\mathbb Z})\, .
\end{equation}

We observed earlier that $H^3({\mathcal U},\, {\mathbb
Z})$ is a pure Hodge structure of weight three.
Therefore, for any smooth compactification
$\overline{\mathcal U}$ of $\mathcal U$, the
homomorphism of mixed Hodge structures
\[
H^3(\overline{\mathcal U},\, {\mathbb Q})\,\longrightarrow
\, H^3({\mathcal U},\, {\mathbb Q})\, ,
\]
induced by the inclusion map ${\mathcal U}\, \hookrightarrow\,
\overline{\mathcal U}$, is surjective \cite[p. 39, Corollaire
(3.2.17)]{De2}. Consequently,
choosing $\overline{\mathcal U}$ to be a smooth
compactification of ${\mathcal M}_X$, from the surjectivity
of the composition homomorphism of mixed Hodge structures
\[
H^3(\overline{\mathcal U},\, {\mathbb Q})\, \longrightarrow\,
H^3({\mathcal M}_X,\, {\mathbb Q})
\,\stackrel{\psi}{\longrightarrow}\, H^3(\mathcal U,\, {\mathbb Q})
\]
we conclude that the homomorphism
\begin{equation}\label{ex.m.3}
\psi_{\mathbb Q}\, :\,
H^3({\mathcal M}_X,\, {\mathbb Q})
\longrightarrow\, H^3(\mathcal U,\, {\mathbb Q})
\end{equation}
induced by the map $\iota$ in Eqn. \eqref{diag.} is surjective.

We noted earlier that $H^4({\mathcal M}_X\, ,
{\mathcal U}\, ,{\mathbb Z})$ is torsionfree. Therefore, from
the surjectivity of $\psi_{\mathbb Q}$ in Eqn. \eqref{ex.m.3}
we conclude that the homomorphism
$\psi$ in Eqn. \eqref{ex.m.2} is surjective.
Consequently, the mixed Hodge structure $H^3({\mathcal M}_X,\,
{\mathbb Z})$ is isomorphic to the pure Hodge structure
$H^3({\mathcal N}_X,\, {\mathbb Z})$ of weight three.
This completes the proof of the proposition.
\end{proof}

Let
\begin{equation}\label{J}
J^2({\mathcal M}_X)\, :=\, H^3({\mathcal M}_X,\, {\mathbb C})/
(F^2H^3({\mathcal M}_X,\, {\mathbb C})+
H^3({\mathcal M}_X,\, {\mathbb Z}))
\end{equation}
be the intermediate Jacobian of the mixed Hodge
structure $H^3({\mathcal M}_X)$ (see \cite[p. 110]{Ca}).
The intermediate Jacobian of any mixed Hodge structure
is a generalized torus \cite[p. 111]{Ca}. Let
\[
J^2({\mathcal N}_X)\, :=\, H^3({\mathcal N}_X,\, {\mathbb C})/
(F^2H^3({\mathcal N}_X,\, {\mathbb C})+
H^3({\mathcal N}_X,\, {\mathbb Z}))
\]
be the intermediate Jacobian for $H^3({\mathcal N}_X,\,
{\mathbb Z})$, which is a complex torus.

\begin{proposition}\label{prop3}
The intermediate Jacobian $J^2({\mathcal M}_X)$ is isomorphic
to $J^2({\mathcal N}_X)$, which is isomorphic to
the Jacobian ${\rm Pic}^0(X)$ of the Riemann surface $X$.
\end{proposition}

\begin{proof}
{}From Proposition \ref{prop2} it follows immediately that
$J^2({\mathcal M}_X)$ is isomorphic to the intermediate Jacobian
$J^2({\mathcal N}_X)$.

On the other hand,
there is a natural isomorphism of $J^2({\mathcal N}_X)$
with $\text{Pic}^0(X)$ \cite[p. 392, Theorem 3]{NR};
this was proved earlier in \cite{MN} for $r\, =\, 2$
(see \cite[p. 1201, Theorem]{MN}). This completes the proof
of the proposition.
\end{proof}

In the next section we will construct a natural polarization
on $J^2({\mathcal M}_X)$.

\section{Polarization on the intermediate Jacobian}

The aim in this section is to show that the Riemann surface $X$
and the variety ${\mathcal M}_X$ determine each other uniquely
up to isomorphisms. We start with a lemma.

\begin{lemma}\label{lem2}
$H^2({\mathcal M}_X,\, {\mathbb Z})\, =\, {\mathbb Z}$.
\end{lemma}

\begin{proof}
Consider the long exact sequence of relative cohomologies
\[
H^2({\mathcal M}_X\, , {\mathcal U}\, ,{\mathbb Z})\,
\longrightarrow\, H^2({\mathcal M}_X\, ,{\mathbb Z})\,
\longrightarrow\, H^2({\mathcal U}\, ,{\mathbb Z})\,
\longrightarrow\,
H^3({\mathcal M}_X\, , {\mathcal U}\, ,{\mathbb Z})
\]
given by $\iota$ in Eqn. \eqref{diag.}. Using
Eqn. \eqref{co.va.}, from this exact sequence
we conclude that $H^2({\mathcal M}_X\, ,{\mathbb Z})\,
=\, H^2({\mathcal U}\, ,{\mathbb Z})$. Therefore, from
the isomorphism in Eqn. \eqref{is.} we have
\begin{equation}\label{2iso.}
H^2({\mathcal M}_X,\, {\mathbb Z})\, =\,
H^2({\mathcal N}_X,\, {\mathbb Z})\, .
\end{equation}
Finally, $H^2({\mathcal N}_X,\, {\mathbb Z})\, =\, {\mathbb Z}$
\cite[p. 582, Proposition 9.13]{AB}. This completes the proof
of the lemma.
\end{proof}

Set $m\, :=\, (r^2-1)(g-1)-3$. Let
\begin{equation}\label{hom.}
F\, :\, (\bigwedge\nolimits^2
H^3({\mathcal M}_X,\, {\mathbb Q}))\otimes
H^2({\mathcal M}_X,\, {\mathbb Q})^{\otimes m}\, \longrightarrow\,
H^{2(r^2-1)(g-1)}({\mathcal M}_X,\, {\mathbb Q})
\end{equation}
be the homomorphism
defined by $(\alpha\bigwedge\beta)\bigotimes \gamma^{\otimes m}
\,\longmapsto \, \alpha\cup\beta\cup \gamma^{\otimes m}$.

\begin{proposition}\label{prop4}
The dimension of the image of the homomorphism $F$, defined
in Eqn. \eqref{hom.}, is one.
\end{proposition}

\begin{proof}
We will use the properties of a
certain moduli space of Higgs bundles over $X$ which is
naturally diffeomorphic to ${\mathcal M}_X$ and is known
as the \textit{Dolbeault moduli space} (see \cite{Si1},
\cite{Ha}).

Let ${\mathcal H}_X$ denote the moduli space parametrizing
all stable Higgs bundles $(E\, ,\theta)$ over $X$ of the
following form:
\begin{itemize}
\item $E$ is a holomorphic vector bundle of rank $r$
with $\bigwedge^r E\, \cong\, {\mathcal O}_X(dx_0)$, and

\item $\text{trace}(\theta) \, \in\, H^0(X,\, K_X)$ vanishes
identically.
\end{itemize}

It is known that the moduli space ${\mathcal H}_X$ is
naturally diffeomorphic to ${\mathcal M}_X$ (see \cite{Ha}). A
quick way to construct the diffeomorphism is the following.

Fix a ramified Galois covering $\phi$ as in Eqn. \eqref{g.co.}.
For any Higgs bundle
$(E\, ,\theta)\, \in\, {\mathcal H}_X$, consider the pulled
back Higgs bundle $((\phi^*E)\bigotimes_{{\mathcal O}_Y}
{\mathcal O}_Y(-dy_0) \, ,\phi^*\theta)$ over $Y$,
where $y_0$ is the point in Eqn. \eqref{pt.}.
This Higgs bundle is polystable of degree zero, because the
pull back of any polystable Higgs bundle by a
Galois covering map remains polystable. Let
$(W\, ,\nabla)$ be the flat vector bundle over $Y$
corresponding to $((\phi^*E)\bigotimes_{{\mathcal O}_Y}
{\mathcal O}_Y(-dy_0) \, ,\phi^*\theta)$ \cite{Hi1},
\cite{Si1}. The flat connection $\nabla$ on $W$ and
the de Rham logarithmic connection on the line bundle
${\mathcal O}_Y(dy_0)$ together induce a logarithmic connection
$\nabla'$ on the vector bundle $W\bigotimes_{{\mathcal O}_Y}{\mathcal
O}_Y(dy_0)$. Let
\[
V\, :=\, (\phi_* (W\otimes_{{\mathcal O}_Y}{\mathcal
O}_Y(dy_0)))^{\text{Gal}(\phi)}
\]
be the invariant direct image on $X$
for the natural action of the
Galois group $\text{Gal}(\phi)$ on the direct image. The logarithmic
connection $\nabla'$ being invariant under the action of the
Galois group, induces a logarithmic connection $\nabla^V$ on $V$.
It is easy to see that $(V\, ,\nabla^V)\,\in\, {\mathcal M}_X$.
Sending any $(E\, ,\theta)$ to $(V\, ,\nabla^V)$, a diffeomorphism
of ${\mathcal H}_X$ with ${\mathcal M}_X$ is obtained. This
diffeomorphism does not depend on the choice of the covering $\phi$.

Since ${\mathcal H}_X$ and ${\mathcal M}_X$ are
diffeomorphic, from Lemma \ref{lem2} we have
$H^2({\mathcal H}_X,\, {\mathbb Q})\, =\, {\mathbb Q}$. Let
\begin{equation}\label{Ga}
\Gamma\, :\, (\bigwedge\nolimits^2 H^3({\mathcal H}_X,
\, {\mathbb Q}))\otimes
H^2({\mathcal H}_X,\, {\mathbb Q})^{\otimes ((r^2-1)(g-1)-3)}\,
\longrightarrow\, H^{2(r^2-1)(g-1)}({\mathcal H}_X,\, {\mathbb Q})
\end{equation}
be the homomorphism
defined by $(\alpha\bigwedge \beta)\bigotimes\gamma^{\otimes m}
\,\longmapsto \, \alpha\cup\beta\cup \gamma^{\otimes m}$.

Comparing the above homomorphism $\Gamma$ with $F$ defined in
Eqn. \eqref{hom.} we conclude that the following lemma
implies that $\dim \text{Image}(F)\, \leq\, 1$.

\begin{lemma}\label{lem3}
The dimension of the image of the homomorphism $\Gamma$, defined
in Eqn. \eqref{Ga}, is at most one.
\end{lemma}

\begin{proof}
To prove this lemma, we consider the \textit{Hitchin map}
\begin{equation}\label{H}
H\, :\, {\mathcal H}_X\, \longrightarrow\, \bigoplus_{i=2}^{r}
H^0(X,\, K^{\otimes i}_X)
\end{equation}
defined by $(E\, ,\theta)\, \longmapsto\, \sum_{i=2}^r
\text{trace}(\theta^i)$ \cite{Hi1}, \cite{Hi2}. This map $H$ is
algebraic and proper \cite{Hi2}, \cite[p. 291, Theorem 6.1]{Ni1}.
The fiber of $H$ over $(0\, , \cdots
\, ,0)$ is known as the \textit{nilpotent cone}. The nilpotent
cone is a finite union of complete subvarieties
of ${\mathcal H}_X$ of complex dimension $(r^2-1)(g-1)$.

Each component of the nilpotent cone defines an element of
$H_{2(r^2-1)(g-1)}({\mathcal H}_X,\, {\mathbb Q})$, and these
elements together generate $H_{2(r^2-1)(g-1)}({\mathcal H}_X,\,
{\mathbb Q})$. This can be proved as follows.

The moduli space ${\mathcal H}_X$ has a natural
structure of a noncompact K\"ahler manifold. It is equipped
with the following holomorphic action of ${\mathbb C}^*$:
\[
\lambda \cdot(E\, ,\theta) \,=\, (E\, ,\lambda\cdot\theta)\, ,
\]
where $\lambda\, \in\, {\mathbb C}^*$ and $(E\, ,\theta)\, \in\,
{\mathcal H}_X$. Also, ${\mathbb C}^*$ acts on
$\bigoplus_{i=2}^{r} H^0(X,\, K^{\otimes i}_X)$ as
\[
\lambda \cdot \sum_{i=2}^r \omega_i\, =\,\sum_{i=2}^r \lambda^i
\cdot\omega_i\, ,
\]
where $\omega_i\, \in\, H^0(X,\, K^{\otimes i}_X)$. The Hitchin
map $H$ is equivariant with respect to these actions of ${\mathbb
C}^*$.

The action of the subgroup $U(1)\, \subset\, {\mathbb C}^*$ on
${\mathcal H}_X$ preserves the K\"ahler metric, and furthermore,
this action is Hamiltonian. The moment map $\mu$ for this action
is defined as follows:
\[
 \mu(E\, , \theta)\,=\,-\frac{1}{2} ||\theta||^2\,=\,-\sqrt{-1}
\int_X \text{trace} (\theta\bigwedge\theta^*)\, ,
\]
where the Hermitian metric on $E$ is the one that satisfies the
Hermitian--Yang--Mills equation. Using Uhlenbeck's compactness
theorem, Hitchin proved that the moment map $\mu$ is proper
\cite{Hi1}. Clearly
$\mu$ is bounded above by $0$. By
a result due to Frankel \cite{Fr}, a proper moment map for a
Hamiltonian action is a perfect Bott--Morse function (the result
of \cite{Fr} is stated for compact K\"ahler manifolds, but it
works in the noncompact case as well). The critical submanifolds
of $\mu$ are the fixed points of the ${\mathbb C}^*$--action on
${\mathcal H}_X$. Using the fact that the Hitchin map $H$ is
equivariant for the actions of ${\mathbb C}^*$ it follows
immediately that
the fixed point set is contained in the nilpotent cone.

We will index the critical submanifolds by ${\mathcal C}_\eta$,
$\eta\, \in\, I$. For each ${\mathcal C}_\eta$,
let $D_\eta$ be
the associated downward Morse flow submanifold, i.e., the
collection of points $x\in {\mathcal H}_X$ whose flow under $-\mu$
converges to ${\mathcal C}_\eta$ (using the K\"ahler form,
any one--form gives a vector field). We will show that
\begin{equation}\label{n.c.i}
 \bigcup_{\eta\in I} D_\eta\,=\, H^{-1}(0,\ldots, 0)\, ,
\end{equation}
which is the nilpotent cone. This result first appeared in
\cite[Theorem 5.2]{Ha-extra}. To prove the equality in Eqn.
\eqref{n.c.i}, first note that given any $z\, \in\, {\mathcal H}_X$
with $H(x_0)\,\neq\, (0,\ldots,0)$, the flow of $z$ diverges as the
flow of $H(x_0)$ diverges (recall that the map $H$ is ${\mathbb
C}^*$--equivariant). On the other hand, the nilpotent cone is
compact, and it is preserved by the action of ${\mathbb C}^*$.
Therefore, for any $z\, \in\, {\mathcal H}_X$ in the nilpotent cone,
the flow of $z$ converges. This proves the equality in Eqn.
\eqref{n.c.i}.

Since the nilpotent cone has finitely many irreducible
components, the index set $I$ is finite, i.e.,
the number of critical submanifolds is finite.
Hence there is a real
number $c\,<\,0$ such that there are no critical points
$(E\, ,\theta)$ with $\mu((E\, ,\theta))\,<\, c$.
Consequently, ${\mathcal H}_X$ retracts to
$\mu^{-1}([c\, ,0])$, and hence it also retracts to the union
$\bigcup_{\eta\in I} D_\eta$, which is the nilpotent cone.

The nilpotent cone is a Lagrangian subvariety of ${\mathcal H}_X$
\cite[p. 648, Th{\'e}or{\`e}me (0.3)]{La}.
Hence each component of the nilpotent cone is a
complete subvariety of (complex) dimension $(r^2-1)(g-1)$.
Therefore, each component of the nilpotent cone defines an element
of $H_{2(r^2-1)(g-1)}({\mathcal H}_X,\, {\mathbb Q})$, and these
elements together generate $H_{2(r^2-1)(g-1)}({\mathcal H}_X,\,
{\mathbb Q})$.

As the next step in the proof of the lemma, we now consider
the moduli space ${\mathcal M}^1_g$ parametrizing all isomorphism
classes of one--pointed compact Riemann
surfaces $(Y\, , y)$ of genus $g$ with $\text{Aut}(Y)
\, =\, e$. This moduli space
is a smooth irreducible quasiprojective variety
of dimension $3g-2$ defined over $\mathbb C$. There is a
universal family of one--pointed Riemann surfaces
\begin{equation}\label{uc}
p\, :\, {\mathcal C}^1_g \, \longrightarrow\, {\mathcal M}^1_g\, .
\end{equation}
The marked points are given by a section of the projection $p$.

Let
\begin{equation}\label{uh}
P\, :\, \widetilde{\mathcal H}\, \longrightarrow\,
{\mathcal M}^1_g
\end{equation}
be the family of moduli spaces of stable Higgs bundles
corresponding to the family of Riemann
surfaces in Eqn. \eqref{uc}. Therefore,
for any one pointed Riemann surface $\underline{x}\, :=\,
(X\, , x_0)\, \in\, {\mathcal M}^1_g$, the fiber
$P^{-1}(\underline{x})$ is the moduli space
${\mathcal H}_X$ parametrizing all stable Higgs bundles
$(E\, ,\theta)$ of rank $r$ with $\bigwedge^r E\, \cong\,
{\mathcal O}_X(dx_0)$ and $\text{trace}(\theta) \, =\, 0$.

For any $j\,\geq\, 1$, consider the direct image on
${\mathcal M}^1_g$
\[
{\mathcal V}_j\, :=\, p_*K^{\otimes i}_{\text{rel}}\, ,
\]
where $K_{\text{rel}}$ is the relative cotangent bundle
on ${\mathcal C}^1_g$ for the projection $p$ in Eqn. \eqref{uc}.
We have the relative Hitchin map
\begin{equation}\label{re.n.co.}
\widetilde{H}\, :\, \widetilde{\mathcal H}\, \longrightarrow\,
\bigoplus_{j=2}^r {\mathcal V}_j
\end{equation}
defined by $(E\, ,\theta)\, \longmapsto\, \sum_{i=2}^r
\text{trace}(\theta^i)$ (see Eqn. \eqref{H}). This
$\widetilde{H}$ is clearly an algebraic morphism.
The relative nilpotent cone is the inverse image
$\widetilde{H}^{-1}(0_{{\mathcal M}^1_g})$, where
$0_{{\mathcal M}^1_g}$ is the zero section of the vector
bundle in Eqn. \eqref{re.n.co.}.

Consider the local system $R^{2(r^2-1)(g-1)}P_*
\underline{{\mathbb C}}$
on ${\mathcal M}^1_g$, where $\underline{{\mathbb C}}$ is the
constant local system on $\widetilde{\mathcal H}$
(constructed in Eqn. \eqref{uh}) with stalk
$\mathbb C$, and $P$ is the projection in
Eqn. \eqref{uh}. Take any $\underline{x}\, =\,(X\,,x)\,
in\, {\mathcal M}^1_g$ and any component $Z$ of the
nilpotent cone of ${\mathcal H}_X$. The element in
\[
H_{2(r^2-1)(g-1)}({\mathcal H}_X,\, {\mathbb Q})\,=\,
H^{2(r^2-1)(g-1)}({\mathcal H}_X,\, {\mathbb Q})^*
\]
given by $Z$ extends uniquely
as a section of $(R^{2(r^2-1)(g-1)}P_*
\underline{{\mathbb C}})^*$, over any contractible analytic
open subset $U_{\underline{x}}\, \subset\, {\mathcal M}^1_g$
containing the point $\underline{x}$. Using the above
construction of
the relative nilpotent cone we conclude that this section
satisfies the following condition:
for each point $\underline{y}\, =\,(Y\,,y)\, in\, U_{\underline{x}}$,
the evaluation of the section
at $\underline{y}$ corresponds to a component of the nilpotent cone
of ${\mathcal H}_Y$.

We noted earlier that
each component of the nilpotent cone in ${\mathcal H}_X$
(there are finitely many of then) defines an
element of $H_{2(r^2-1)(g-1)}({\mathcal H}_X,\, {\mathbb Q})$,
and these elements together generate
$H_{2(r^2-1)(g-1)}({\mathcal H}_X,\, {\mathbb Q})$. This
together with the above observation on
relative nilpotent cone
imply that the monodromy of the local system
$R^{2(r^2-1)(g-1)}P_*{\underline{\mathbb C}}$ is a quotient group
of the permutation group of the components of the nilpotent cone.
In particular, the monodromy of the local system
$R^{2(r^2-1)(g-1)}P_*\underline{{\mathbb C}}$ is a finite group.

Fix a base point $\underline{x}_0\, :=\, (X_0\, ,x_0)\, \in\,
{\mathcal M}^1_g$. Let $G_{\mathbb Z}$
(respectively, $G_{\mathbb C}$)
denote the group of all automorphisms of $H^1(X_0,\, {\mathbb Z})$
(respectively, $H^1(X_0,\, {\mathbb C})$) preserving
the symplectic pairing given by the cup product. Choosing
a symplectic basis of $H^1(X_0,\, {\mathbb Z})$, the groups
$G_{\mathbb Z}$ and $G_{\mathbb C}$ get identified
with $\text{Sp}(2g, {\mathbb Z})$ and $\text{Sp}(2g, {\mathbb C})$
respectively.

Consider the local system
$R^{1}p_*\underline{\underline{\mathbb Z}}$ on
${\mathcal M}^1_g$, where $p$ is the projection in Eqn. \eqref{uc},
and $\underline{\underline{\mathbb Z}}$ is the constant
local system on ${\mathcal C}^1_g$ with stalk $\mathbb Z$.
Using its monodromy, the group $G_{\mathbb Z}$ is a quotient
of the fundamental group
\[
\Gamma^1_g \, :=\, \pi_1({\mathcal M}^1_g,\, \underline{x}_0)\, .
\]
This group $\Gamma^1_g$
is known as the \textit{mapping class group}, and the kernel of
the projection of $\Gamma^1_g$ to $G_{\mathbb Z}$ is
known as the \textit{Torelli group}.

The cohomology algebra $H^*({\mathcal H}_X,\, {\mathbb C})$
is generated by the K\"unneth components of a universal
vector bundle over $X\times {\mathcal H}_X$
\cite[p. 73, Theorem 7]{Mar}; this was proved earlier in
\cite{HT} for $r\, =\, 2$ (see \cite[p. 641, (6.1)]{HT}). From
this it follows immediately that for any $i\, \geq\, 0$,
the local system $R^{i}P_*\underline{{\mathbb C}}$ is a quotient
of some local system on ${\mathcal M}^1_g$ of the form
\[
{\mathcal W}\, :=\, \bigoplus_{j=1}^\ell\big(\big(
\bigoplus_{i=0}^2 R^{i}p_*\underline{\underline{\mathbb C}}
\big)^{\oplus a_j}\big)^{\otimes b_j}\, ,
\]
where $a_j,b_j\, \in\, {\mathbb N}$, and
$p$ is the projection in Eqn. \eqref{uc}. We note that using
the marked point, a universal Higgs bundle can be rigidified.
Since both
$R^{0}p_*\underline{\underline{\mathbb C}}$ and
$R^{2}p_*\underline{\underline{\mathbb C}}$ are constant local
systems, and $\text{Sp}(2g, {\mathbb Z})$ is Zariski
dense in $\text{Sp}(2g, {\mathbb C})$ \cite[p. 179]{Bo},
we conclude the following:
\begin{enumerate}
\item{} For all $i\, \geq\, 0$,
the Torelli group is in the kernel of the monodromy
representation
\begin{equation}\label{m.h.}
\Gamma^1_g\, \longrightarrow\, \text{Aut}(
(R^{i}P_*\underline{{\mathbb C}})_{\underline{x}_0})
\end{equation}
of the mapping class group for the local
system $R^{i}P_*\underline{{\mathbb C}}$,
where $\underline{x}_0\, \in\, {\mathcal M}^1_g$ is the base
point. Therefore, the homomorphism in Eqn. \eqref{m.h.}
factors through the quotient $G_{\mathbb Z}$ of $\Gamma^1_g$.

\item{} The homomorphism
\begin{equation}\label{m.h2.}
G_{\mathbb Z}\,\longrightarrow
\,\text{Aut}((R^{i}P_*\underline{{\mathbb C}})_{\underline{x}_0})
\end{equation}
obtained from Eqn. \eqref{m.h.} extends to a representation of
$G_{\mathbb C}$.
\end{enumerate}
To prove the second assertion,
note that as $G_{\mathbb Z}$ is Zariski dense in
$G_{\mathbb C}$, the kernel of the surjective homomorphism of
$G_{\mathbb Z}$--modules
\[
{\mathcal W}_{\underline{x}_0}\, \longrightarrow\,
(R^{i}P_*\underline{{\mathbb C}})_{\underline{x}_0}
\]
is preserved by the action of $G_{\mathbb C}$ on
${\mathcal W}_{\underline{x}_0}$, thus inducing and action
of $G_{\mathbb C}$ on the quotient.

We noted earlier that the monodromy of the local system
$R^{2(r^2-1)(g-1)}P_*\underline{\mathbb C}$ is a finite group.
In other words, the monodromy representation for
$R^{2(r^2-1)(g-1)}P_*\underline{\mathbb C}$ factors through
a finite quotient of $G_{\mathbb Z}\,\cong\,
\text{Sp}(2g, {\mathbb Z})$. Any finite index
subgroup of $\text{Sp}(2g, {\mathbb Z})$ is Zariski
dense in $\text{Sp}(2g, {\mathbb C})$ \cite[p. 179]{Bo}.
Since a Zariski dense subgroup is in the kernel of the
homomorphism
\[
G_{\mathbb C}\, \longrightarrow\, \text{Aut}((R^{2(r^2-1)(g-1)}
P_*\underline{{\mathbb C}})_{\underline{x}_0})
\]
obtained by extending the monodromy representation of
$G_{\mathbb Z}$ in Eqn. \eqref{m.h2.}
for the local system $R^{2(r^2-1)(g-1)}
P_*\underline{{\mathbb C}}$, we conclude that the above
homomorphism is the trivial homomorphism.
In other words, the monodromy of $R^{2(r^2-1)(g-1)}
P_*\underline{{\mathbb C}}$ is trivial.
Consequently, the local
system $R^{2(r^2-1)(g-1)}P_*\underline{{\mathbb C}}$ is
a constant one.

Consider the homomorphism $\Gamma$ constructed in Eqn. \eqref{Ga}
with coefficients in $\mathbb C$ instead of $\mathbb Q$.
The pointwise construction of it (over the points of
${\mathcal M}^1_g$) yields a homomorphism of local systems
\begin{equation}\label{tga}
\widetilde{\Gamma}\, :\,
(\bigwedge\nolimits^2 R^{3}P_*\underline{\mathbb C}) \otimes
(R^{2}P_*\underline{\mathbb C})^{\otimes ((r^2-1)(g-1)-3)}\,
\longrightarrow\, R^{2(r^2-1)(g-1)}P_*\underline{\mathbb C}\, .
\end{equation}

{}From Lemma \ref{lem2} it follows that
$R^{2}P_*\underline{\mathbb C}$ is a constant local
system of rank one.

There is a canonical isomorphism
$H^3({\mathcal N}_X,\, {\mathbb Z})\, =\,
H^1(X,\, {\mathbb Z})$, where ${\mathcal N}_X$ is the moduli
space of stable vector bundles over $X$
defined in Section \ref{sec3} \cite[p. 392, Theorem 3]{NR}.
This isomorphism is constructed using the K\"unneth component
of the second Chern class of a universal
vector bundle over $X\times {\mathcal N}_X$.
Therefore, using Proposition \ref{prop2} it follows that
\[
R^{3}P_*\underline{\mathbb C}\, =\, R^1p_*
\underline{\underline{\mathbb C}}\, ,
\]
where $p$ is the projection in Eqn. \eqref{uc}.

If we identify $G_{\mathbb Z}$ with $\text{Sp}(2g, {\mathbb Z})$
by choosing a symplectic basis of $H^1(X_0,\, {\mathbb Z})$,
then $R^1p_*\underline{\underline{\mathbb C}}$ gets identified
with the local system on ${\mathcal M}^1_g$ associated to the
standard representation of $\text{Sp}(2g,{\mathbb Z})$. From
this it follows that the local system $\bigwedge^2
R^1p_*\underline{\underline{\mathbb C}}$ decomposes as
\[
\bigwedge\nolimits^2 R^1p_*\underline{\underline{\mathbb C}}
\, =\,{\mathbb L}^0\oplus {\mathbb L}^1\, ,
\]
where ${\mathbb L}^0$ is a constant local system of rank one
and ${\mathbb L}^1$ is an irreducible local system of rank
$g(2g-1)-1$ on ${\mathcal M}^1_g$. This decomposition corresponds
to the decomposition of the $\text{Sp}(2g,{\mathbb C})$--module
$\bigwedge^2 {\mathbb C}^{2g}$ as a direct sum of the
trivial $\text{Sp}(2g,{\mathbb C})$--module of dimension one
with an irreducible $\text{Sp}(2g,{\mathbb C})$--module.
Since ${\mathbb L}^1$ is irreducible of rank
at least two, and
both $R^{2(r^2-1)(g-1)}P_*\underline{\mathbb C}$ and
$R^{2}P_*\underline{\mathbb C}$ are constant local systems,
the restriction of the homomorphism $\widetilde{\Gamma}$
(constructed in Eqn. \eqref{tga}) to the sublocal system
\[
{\mathbb L}^1 \otimes
(R^{2}P_*\underline{\mathbb C})^{\otimes ((r^2-1)(g-1)-3)}\,
\subset\,(\bigwedge\nolimits^2 R^{3}P_*\underline{\mathbb C})
\otimes
(R^{2}P_*\underline{\mathbb C})^{\otimes ((r^2-1)(g-1)-3)}
\]
vanishes identically.

Therefore, the homomorphism $\widetilde{\Gamma}$ factors
through the one dimensional quotient local system
${\mathbb L}^0\bigotimes
(R^{2}P_*\underline{\mathbb C})^{\otimes ((r^2-1)(g-1)-3)}$.
This immediately implies that the dimension of the
image of the homomorphism
$\Gamma$ defined in Eqn. \eqref{Ga} is at most one.
This completes the proof of the lemma.
\end{proof}

Continuing with the proof of Proposition \ref{prop4}, from
Lemma \ref{lem3} it follows that
\[
\dim \text{Image}(F)\, \leq\, 1\, .
\]
We will complete the proof of the proposition by showing that
$\text{Image}(F)\, \not=\, 0$.

Let
\begin{equation}\label{g1}
\gamma\, \in\, H^2({\mathcal N}_X,\, {\mathbb Q})
\end{equation}
be an ample class of the smooth projective variety
${\mathcal N}_X$ (defined in Section \ref{sec3}). The
Hard Lefschetz theorem for ${\mathcal N}_X$ says that
there are
\begin{equation}\label{ab}
\alpha\, , \beta\, \in\, H^3({\mathcal N}_X,\, {\mathbb Q})
\end{equation}
such that
\begin{equation}\label{nz}
0\, \not=\, \alpha\cup\beta\cup
\gamma^{\otimes ((r^2-1)(g-1)-3)}
\,\in\, H^{2(r^2-1)(g-1)}({\mathcal N}_X,\, {\mathbb Q})
\,=\, {\mathbb Q}
\end{equation}
(recall that $\dim_{\mathbb C}{\mathcal N}_X\, =\,
(r^2-1)(g-1)$).

Let
\begin{equation}\label{no1}
\widetilde{\alpha}\, , \widetilde{\beta}\, \in\,
H^3({\mathcal M}_X,\, {\mathbb Q})
\end{equation}
be the elements given by $\alpha$ and $\beta$
respectively (in Eqn. \eqref{ab}) using the isomorphism
$H^3({\mathcal N}_X,\, {\mathbb Q})\, =\,
H^3({\mathcal M}_X,\, {\mathbb Q})$ in Proposition
\ref{prop2}. Similarly, let
\begin{equation}\label{no2}
\widetilde{\gamma}\, \in\,
H^2({\mathcal M}_X,\, {\mathbb Q})
\end{equation}
be the element given by $\gamma$ using the isomorphism
in Eqn. \eqref{2iso.}.

We will show that
\begin{equation}\label{idc}
F((\widetilde{\alpha}\wedge\widetilde{\beta})\otimes
\widetilde{\gamma}^{\otimes ((r^2-1)(g-1)-3)})
\, \not=\, 0\, ,
\end{equation}
where $F$ is the homomorphism in Eqn. \eqref{hom.}.

The projection $\Phi$ in Eqn. \eqref{Ph} has a natural
$C^\infty$ section given the unique unitary
logarithmic connection singular over $x_0$
on any stable vector bundle. We will briefly explain the
construction of this section. Fix a ramified Galois covering
$\phi$ as in Eqn. \eqref{g.co.}.
Given any stable vector
bundle $E\, \in\, {\mathcal N}_X$, consider the vector
bundle $(\phi^*E)\bigotimes_{{\mathcal O}_Y}
{\mathcal O}_Y(-dy_0)$ over $Y$, where $y_0$ is the
point in Eqn. \eqref{pt.}. This vector bundle is
polystable of degree zero. Hence it admits a unique unitary
flat connection $\nabla$ \cite{NS}. Let $\nabla'$ be the
logarithmic connection on
\[
\phi^*E\, =\, ((\phi^*E)\otimes_{{\mathcal O}_Y}
{\mathcal O}_Y(-dy_0))\otimes_{{\mathcal O}_Y}{\mathcal O}_Y(dy_0)
\]
induced by $\nabla$ and the de Rham logarithmic connection
on ${\mathcal O}_Y(dy_0)$. The above mentioned unique unitary
connection on $E$ is the
unique logarithmic connection $\nabla''$ on $E$ such that the
pullback $\phi^*\nabla''$ coincides with $\nabla'$. It is
easy to see that $(E\, ,\nabla'')\, \in\, {\mathcal M}_X$.

Let
\[
\tau\, :\, {\mathcal N}_X\, \longrightarrow\, {\mathcal U}
\]
be the smooth section of $\Phi$ that associates to any stable
vector bundle $E$ the unique unitary logarithmic
connection $\nabla''$ on it. Let
\[
\widetilde{\tau}\, \in\, H_{2(r^2-1)(g-1)}({\mathcal M}_X,
\, {\mathbb Q})
\]
be the homology class defined by the image of $\tau$ using
the inclusion of ${\mathcal U}$ in ${\mathcal M}_X$.

Consider the natural duality pairing between the
cohomology $H^{2(r^2-1)(g-1)}(
{\mathcal M}_X, \, {\mathbb Q})$ and the homology
$H_{2(r^2-1)(g-1)}({\mathcal M}_X, \, {\mathbb Q})$.
It is straight--forward to check that
\begin{equation}\label{id.p.}
F((\widetilde{\alpha}\wedge\widetilde{\beta})\otimes
\widetilde{\gamma}^{\otimes ((r^2-1)(g-1)-3)})(\widetilde{\tau})
\,=\, (\alpha\cup\beta\cup \gamma^{\otimes ((r^2-1)(g-1)-3)})
\cap [{\mathcal N}_X]\, \in\, {\mathbb Q}\, ,
\end{equation}
where $\alpha$, $\beta$ and $\gamma$ are the cohomology
classes in Eqn. \eqref{ab} and Eqn. \eqref{g1} respectively,
and $\widetilde{\alpha}$, $\widetilde{\beta}$
and $\widetilde{\alpha}$ are constructed in Eqn. \eqref{no1}
and Eqn. \eqref{no2} respectively. Using Eqn. \eqref{id.p.}
we conclude that the assertion in Eqn. \eqref{idc} follows from
Eqn. \eqref{nz}. Consequently, $\text{Image}(F)\, \not=\, 0$. We
have already shown that $\dim \text{Image}(F)\, \leq\, 1$.
This completes the proof of the proposition.
\end{proof}

Now we are in a position to prove our main result.

\begin{theorem}\label{thm1}
Let $(X\, ,x_0)$ and $(Y\, ,y_0)$ be two compact
connected one--pointed Riemann
surfaces of genus $g$, with $g\, \geq\, 3$. Let ${\mathcal
M}_X$ and ${\mathcal M}_Y$ be the corresponding moduli
spaces of connections defined as in Section \ref{sec1}.
The two varieties ${\mathcal M}_X$ and ${\mathcal M}_Y$
are isomorphic if and only if the two Riemann surfaces
$X$ and $Y$ are isomorphic.
\end{theorem}

\begin{proof}
{}From Proposition \ref{prop4} we know that
$\dim \text{Image}(F)\, =\, 1$.
Consequently, fixing a generator of
$\text{Image}(F)$, and also fixing
a generator of $H^2({\mathcal M}_X,\, {\mathbb Q})$,
the homomorphism $F$ gives a cohomology class
\begin{equation}\label{theta}
\theta\, \in\, \bigwedge\nolimits^2
H^3({\mathcal M}_X,\, {\mathbb Q})^*
\, =\, H^2(J^2({\mathcal M}_X),\, {\mathbb Q})\, .
\end{equation}

Using the diagram in Eqn. \eqref{diag.}, we have an isomorphism
of $J^2({\mathcal M}_X)$ with $J^2({\mathcal N}_X)$
(Proposition \ref{prop3}). Consider the N{\'e}ron--Severi group
\[
\text{NS}(J^2({\mathcal N}_X))\, :=\,
H^{1,1}(J^2({\mathcal N}_X))\cap H^2(J^2({\mathcal N}_X),\,
{\mathbb Z})\, .
\]
There is a natural element
\[
\widetilde{\theta}\,\in\,
\text{NS}(J^2({\mathcal N}_X))_{\mathbb Q} \, :=\,
\text{NS}(J^2({\mathcal N}_X))\otimes_{\mathbb Z} {\mathbb Q}
\,=\,
H^{1,1}(J^2({\mathcal N}_X))\cap H^2(J^2({\mathcal N}_X),\,
{\mathbb Q})
\]
defined by
\[
\alpha\wedge \beta\, \longmapsto\, (\alpha\cup\beta
\cup {\gamma}^{\otimes ((r^2-1)(g-1)-3)})\cap [{\mathcal N}_X]
\, \in\, {\mathbb Q}\, ,
\]
where $\alpha\, ,\beta\, \in\, H^3({\mathcal N}_X,\, {\mathbb Q})$,
and $\gamma$ is the ample generator of
$H^2({\mathcal N}_X,\, {\mathbb Z})$. The
Hard Lefschetz theorem for ${\mathcal N}_X$
says that $\widetilde{\theta}\,\not=\, 0$.
We will show that $\widetilde{\theta}$ is a scalar multiple
of a principal polarization on $J^2({\mathcal N}_X)$.

Let
\[
\overline{p} \,:\, {\mathcal J}\, \longrightarrow\,
{\mathcal M}^1_g
\]
be the universal Jacobian for the universal family of
one--pointed Riemann surfaces in Eqn. \eqref{uc}. It
is known that any section of the local system
$R^2\overline{p}_*{\mathbb Q}$ which is a Hodge cycle
over every point of ${\mathcal M}^1_g$ (i.e., an element
of $\text{NS}(\text{Pic}^0(Y))_{\mathbb Q}$ for every
Riemann surface $Y$), must be a scalar multiple of the
section of $R^2\overline{p}_*{\mathbb Q}$ given by the
theta line bundle on the Jacobians;
this is well--known (see \cite[p. 712]{BiN}).

We have $J^2({\mathcal N}_X)\, =\, {\rm Pic}^0(X)$
(see \cite[p. 392, Theorem 3]{NR}). Therefore,
in view of the above remark, the cohomology class
$\widetilde{\theta}$ is a rational multiple of
a principal polarization on $J^2({\mathcal N}_X)$.

Given a nonzero
rational multiple, say $\theta'$, of a principal polarization
on an abelian variety $A$, there is a unique way to scale
$\theta'$ to get back the principal polarization on $A$. Indeed,
this is an immediately consequence of the following
two properties of a principal polarization:
\begin{itemize}
\item{} a principal polarization on $A$ is an ample class on $A$,
and
\item{} a principal polarization is indivisible as an element
of $H^2(A,\, {\mathbb Z})$.
\end{itemize}
Note that a free $\mathbb Z$--module of rank one has
exactly two generators.

Let $\widetilde{\widetilde{\theta}}$ be the unique principal
polarization on $J^2({\mathcal N}_X)$ which is a rational
multiple of $\widetilde{\theta}$. Therefore,
the principally polarized abelian variety
$(J^2({\mathcal N}_X)\, ,\widetilde{\widetilde{\theta}})$
is isomorphic to ${\rm Pic}^0(X)$ equipped with the canonical
principal polarization on it given by the theta line bundle.

Comparing $\widetilde{\theta}$ with $\theta$
defined in Eqn. \eqref{theta} we conclude that the
above mentioned identification of
$J^2({\mathcal N}_X)$ with $J^2({\mathcal M}_X)$ takes
$\widetilde{\widetilde{\theta}}$ to a nonzero multiple of
$\theta$. Now using the
earlier remark that a principal polarization can be
recovered from any nonzero
multiple of it we conclude that $\theta$
gives a principal polarization $\widehat{\theta}$
on $J^2({\mathcal M}_X)$. Furthermore,
the principally polarized abelian
variety $(J^2({\mathcal M}_X)\, ,\widehat{\theta})$ is isomorphic
to the principally polarized abelian variety
$(J^2({\mathcal N}_X)\, ,\widetilde{\widetilde{\theta}})$.

We have noted earlier that
the principally polarized abelian variety
$(J^2({\mathcal N}_X)\, ,\widetilde{\widetilde{\theta}})$
is isomorphic to ${\rm Pic}^0(X)$ equipped with the canonical
principal polarization on it given by the theta line bundle.
Therefore, from the classical Torelli theorem, which says
that the isomorphism class of the principally polarized
abelian variety ${\rm Pic}^0(X)$ equipped with the
principal polarization given by the theta line bundle
determines the Riemann surface $X$
uniquely up to an isomorphism, we
conclude that the isomorphism class of the variety
${\mathcal M}_X$ determines the Riemann surface $X$
uniquely up to an isomorphism.

We saw in Remark \ref{rem1} that the isomorphism class
of the variety ${\mathcal M}_X$ is independent of the
choice of the base point $x_0$. This implies that the
isomorphism class of the variety ${\mathcal M}_X$ is determined
uniquely by the isomorphism class of the Riemann surface $X$.
Therefore, the proof of the theorem is complete.
\end{proof}

In the next section we will prove a similar result
for the moduli space
of $\text{GL}(n,{\mathbb C})$--connections.

\section{Moduli of $\text{GL}(n,{\mathbb C})$--connections}

As before, fix two mutually coprime integers $r$ and $d$, with
$r\, \geq\, 2$. Let $(X\, ,x_0)$ be a one--pointed compact
Riemann surface of genus $g$, with $g\,\geq\, 3$.

Let $\widehat{\mathcal M}_X$ be the moduli space parametrizing
all logarithmic connections $(E\, ,D)$ over $X$ singular
exactly over $x_0$ and satisfying the following two conditions:
\begin{itemize}
\item $\text{rank}(E)\, =\, r$, and

\item $\text{Res}(D,x_0) \, =\,
-\frac{d}{r}\text{Id}_{E_{x_0}}$.
\end{itemize}
{}From Eqn. \eqref{de-re} it follows that if
$(E\, ,D)\, \in\, \widehat{\mathcal M}_X$, then
$\text{degree}(E)\,=\, d$.
The moduli space $\widehat{\mathcal M}_X$ is an irreducible
smooth quasiprojective variety defined over $\mathbb C$. Its
dimension is $2r^2(g-1)+2$.

A surjective morphism
\[
p_0\, :\, \widehat{\mathcal M}_X\, \longrightarrow\, A_0\, ,
\]
where $A_0$ is an abelian variety defined over
$\mathbb C$, will be called
\textit{universal} if for any morphism
\begin{equation}\label{mor}
p\, :\, \widehat{\mathcal M}_X\, \longrightarrow\, A\, ,
\end{equation}
where $A$ is any abelian variety
defined over $\mathbb C$, there is a unique morphism
\[
\gamma\, :\, A_0\, \longrightarrow\, A
\]
such that $\gamma{\circ}p_0\, =\, p$.

\begin{proposition}\label{propu}
The morphism
\begin{equation}\label{prou}
p_0\, :\, \widehat{\mathcal M}_X\, \longrightarrow\,
{\rm Pic}^d(X)
\end{equation}
defined by $(E\, ,D)\,\longmapsto\, \bigwedge^r E$
is universal.
\end{proposition}

\begin{proof}
{}From \cite[p. 787, Lemma 2.3]{BR} it follows immediately
that the morphism $p_0$ in Eqn. \eqref{prou} is surjective.
To prove that it is universal,
let $A$ be an abelian variety defined over $\mathbb C$
and $p$ a morphism to $A$ as in Eqn. \eqref{mor}. Fix
any closed point $z\,\in\, {\rm Pic}^d(X)$. Let
\begin{equation}\label{pz}
p^z\, :\, p^{-1}_0(z) \, \longrightarrow\, A
\end{equation}
be the restriction of $p$ to the closed subvariety
$p^{-1}_0(z)\, \subset\, \widehat{\mathcal M}_X$.

Fix a holomorphic line bundle
$L_0$ over $X$ of degree zero such that
the line bundle $(L^*_0)^{\otimes r}
\bigotimes_{{\mathcal O}_X}{\mathcal O}_X(dx_0)$
corresponds to the point $z\, \in\, \text{Pic}^d(X)$.
Consider the moduli space ${\mathcal M}_X$
defined in Section \ref{sec1}. Given any
$(E\, ,D)\, \in\, p^{-1}_0(z)$, there is a unique
holomorphic connection $\nabla$ on $L_0$ (the
connection $\nabla$ depends on $D$) such that
\begin{equation}\label{h.i.}
(E\otimes L_0\, ,D\otimes\text{Id}_{L_0}+
\text{Id}_{E}\otimes\nabla)\, \in\, {\mathcal M}_X\, .
\end{equation}
Let
\begin{equation}\label{ga.}
\gamma\, :\, p^{-1}_0(z)\, \longrightarrow\, {\mathcal M}_X
\end{equation}
be the morphism that sends any $(E\, ,D)\, \in\,
p^{-1}_0(z)$ to $(E\otimes L_0\, ,D\bigotimes\text{Id}_{L_0}+
\text{Id}_{E}\bigotimes\nabla)$ in Eqn. \eqref{h.i.}.
This map $\gamma$ is algebraic, and it
makes $p^{-1}_0(z)$ an affine bundle
over ${\mathcal M}_X$ for the trivial vector bundle over
${\mathcal M}_X$ with fiber $H^0(X,\, K_X)$. In other
words, $H^0(X,\, K_X)$ acts freely transitively on
each of the fibers of $\gamma$.

Any algebraic morphism $f\, :\, {\mathbb C}{\mathbb P}^1\,
\longrightarrow\, A'$, where $A'$ is a complex abelian variety,
is a constant morphism. Indeed, the cotangent bundle of $A'$
is generated by global sections, and
\[
H^0({\mathbb C}{\mathbb P}^1,\, K_{{\mathbb C}{\mathbb P}^1})
\,=\, 0\, ,
\]
implying that the differential of $f$ vanishes. Therefore, any
algebraic morphism from a Zariski open subset of
${\mathbb C}{\mathbb P}^1$ to $A'$ is a constant morphism.

{}From the above observation it follows that
the restriction of the morphism $p^z$
(defined in Eqn. \eqref{pz}) to any fiber of
$\gamma$ (defined in Eqn. \eqref{ga.}) is a constant
morphism. Consequently, there is unique morphism
\begin{equation}\label{ga0.}
\gamma_0\, :\, {\mathcal M}_X\,\longrightarrow\, A
\end{equation}
such that $\gamma_0\circ \gamma\, =\, p^z$.

Let
\begin{equation}\label{ga1.}
\gamma_1\, :\, {\mathcal U}\,\longrightarrow\, A
\end{equation}
be the restriction of $\gamma_0$ to the open dense subset
$\mathcal U\, \subset\, {\mathcal M}_X$
defined in Eqn. \eqref{cU}. Since
any fiber of the morphism $\Phi$ in Eqn. \eqref{Ph} is
an affine space, the restriction of $\gamma_1$ to any fiber
of $\Phi$ is a constant morphism. Therefore, there is
a unique morphism
\begin{equation}\label{ga2.}
\gamma_2\, :\, {\mathcal N}_X\,\longrightarrow\, A
\end{equation}
such that $\gamma_2\circ\Phi\, =\, \gamma_1$, where
$\gamma_1$ is defined in Eqn. \eqref{ga1.}.

The variety ${\mathcal N}_X$ is known to be unirational,
in fact, it is rational \cite[p. 520, Theorem 1.2]{KS}.
Hence a nonempty Zariski open subset of
it is covered by the images of ${\mathbb C}{\mathbb P}^1$.
This implies that the morphism $\gamma_2$ in Eqn.
\eqref{ga2.} is a constant one.

Therefore, $\gamma_1$
in Eqn. \eqref{ga1.} is a constant morphism. Since
$\mathcal U$ is Zariski dense in ${\mathcal M}_X$,
the morphism $\gamma_0$ in Eqn. \eqref{ga0.} is a constant
one. Consequently, $p^z$ in Eqn. \eqref{pz} is a constant
morphism. This implies that the morphism $p$ factors
through a morphism $\text{Pic}^d(X)\,\longrightarrow\, A$.
Hence the morphism
$p_0$ in Eqn. \eqref{prou} is universal. This completes
the proof of the proposition.
\end{proof}

{}From Proposition \ref{propu} it follows that
the isomorphism class of the variety $\widehat{\mathcal M}_X$
determines the pair ${\rm Pic}^d(X)$ and the projection
$p_0$ in Eqn. \eqref{prou} up to an automorphism of
${\rm Pic}^d(X)$. Therefore,
the isomorphism class of the variety $\widehat{\mathcal M}_X$
determines the isomorphism class of the variety
$p^{-1}_0(z)$ for some $z\, \in\, \text{Pic}^0(X)$. More
precisely, for any point $w\, \in\, \widehat{\mathcal M}_X$,
the isomorphism class of the variety
$p^{-1}_0(p_0(w))$ is determined by the
isomorphism class of the variety $\widehat{\mathcal M}_X$.

Since the morphism $\gamma$ in Eqn. \eqref{ga.} is an
affine fibration, in particular, the fibers are contractible,
the induced homomorphism
\[
\gamma^*\, :\, H^i({\mathcal M}_X,\, {\mathbb Z})
\, \longrightarrow\, H^i(p^{-1}_0(z),\, {\mathbb Z})
\]
is an isomorphism for all $i\, \geq\, 0$. Therefore,
$J^2({\mathcal M}_X)\, \cong\, J^2(p^{-1}_0(z))$.
Furthermore, the cup product
\[
(\bigwedge\nolimits^2 H^3(p^{-1}_0(z)\, {\mathbb Q}))\otimes
H^2(p^{-1}_0(z),\, {\mathbb Q})^{\otimes ((r^2-1)(g-1)-3)}
\, \longrightarrow\,
H^{2(r^2-1)(g-1)}(p^{-1}_0(z),\, {\mathbb Q})
\]
defined as in Eqn. \eqref{hom.} gives, after rescaling,
a principal polarization on $J^2(p^{-1}_0(z))$. The
resulting principally polarized abelian variety is clearly
isomorphic to $(J^2({\mathcal M}_X)\, ,\widehat{\theta})$
constructed in the proof of Theorem \ref{thm00}.

Therefore, using Theorem \ref{thm1}, we have the
following theorem:

\begin{theorem}\label{thm000}
The isomorphism class of the variety
$\widehat{\mathcal M}_X$ determines the Riemann surface
$X$ uniquely up to an isomorphism.
\end{theorem}

\begin{remark}\label{rem.f}
{\rm The isomorphism class of the variety $\widehat{\mathcal M}_X$
is independent of the choice of the point $x_0$. This can be shown
using the morphism defined in Remark \ref{rem1}. More precisely,
if $\widehat{\mathcal M}^{x_1}_X$ is the moduli space of logarithmic
connection obtained by replacing $x_0$ with a different
point $x_1$ in the
construction of $\widehat{\mathcal M}_X$, then the morphism
$$
\widehat{\mathcal M}_X\, \longrightarrow\, \widehat{\mathcal M}^{x_1}_X
$$
defined by
$$
(E\, ,D)\, \longmapsto\, (E\otimes L\, , D\otimes \text{Id}_{L}+
\text{Id}_E \otimes D_0)
$$
is an algebraic isomorphism of varieties, where
$(L\, , D_0)$ is defined in Remark \ref{rem1}.

Also, the
complex manifold
$\widehat{\mathcal M}_X$ is naturally biholomorphic to the subset of
$\text{Hom}(\pi_1(X\setminus\{x_0\})\, , \text{GL}(r, {\mathbb C}))$
parametrizing all representations that sends the anticlockwise
oriented loop around $x_0$ to
$\exp(2\pi\sqrt{-1}d/r)\cdot I_{r\times r}\,\in\,
\text{GL}(r, {\mathbb C})$. As before, the biholomorphism
is obtained by sending a connection to the corresponding monodromy
representation.}
\end{remark}

\begin{remark}\label{ra1}
{\rm It can be shown that Theorem \ref{thm000} fails for rank
one case. To prove this, let $\widehat{\mathcal M}^1$ denote the
moduli space of rank one regular connections on $X$. So there
is an algebraic morphism
$$
\psi\, :\, \widehat{\mathcal M}^1\, \longrightarrow\,
{\rm Pic}^0(X)\, =:\, T
$$
that sends any $(L\, ,D)\, \in\, \widehat{\mathcal M}^1$
to $L\, \in\, {\rm Pic}^0(X)$.

Let $\Omega^1_T$ denote the holomorphic cotangent bundle
of the abelian variety $T\,=\,{\rm Pic}^0(X)$.
Since any two holomorphic connections on a given holomorphic
line bundle over $X$ differ by a holomorphic one--form on $X$,
the morphism $\psi$ makes $\widehat{\mathcal M}^1$ a
$\Omega^1_T$--torsor. In other words, for each point $t\,
\in\, T$, the fiber $\psi^{-1}(t)$ is an affine space for the
vector space $(\Omega^1_T)_x \, =\, H^0(X,\, K_X)$. Isomorphism
classes of $\Omega^1_T$--torsors are parametrized by
$H^1(T,\, \Omega^1_T)$. It is known that the above
$\Omega^1_T$--torsor $\widehat{\mathcal M}^1$ corresponds to
$\Theta/2$, where $\Theta\, \in\, H^1(T,\, \Omega^1_T)$
is the class of a theta divisor on $T$ \cite[p. 308,
Theorem 2.11]{BR2}.

Elements of $H^1(T,\, \Omega^1_T)$ correspond to the
translation invariant vector bundle homomorphisms from
the smooth tangent bundle $T^{0,1}T$ of type $(0\, ,1)$ to
the cotangent bundle $\Omega^1_T$; here translation invariant
means invariant for the translation action of the torus
$T$ on itself. The class $\Theta/2\, \in\, H^1(T,\, \Omega^1_T)$
corresponds to an isomorphism of $\Omega^1_T$ with
$T^{0,1}T$.

Let ${\rm Aut}(\Omega^1_T)$ denote the
group of all holomorphic automorphisms of the vector bundle
$\Omega^1_T$. Consider the natural action of ${\rm Aut}(\Omega^1_T)$
on the space of all translation invariant vector
bundle homomorphisms from $T^{0,1}T$ to $\Omega^1_T$;
the action of any $\tau\, \in\,
{\rm Aut}(\Omega^1_T)$ sends a homomorphism
$\beta$ to $\tau\circ\beta$. It is easy to see that
the action of ${\rm Aut}(\Omega^1_T)$ on the
space of all translation invariant
isomorphisms of $T^{0,1}T$ with $\Omega^1_T$ is transitive.

This implies that if
$$
\psi_i\, :\, {\mathcal A}_i\, \longrightarrow\, T\, ,
$$
$i\, =\, 1,2$, are two algebraic
$\Omega^1_T$--torsors such that the translation
invariant homomorphisms, from $T^{0,1}T$
to $\Omega^1_T$, corresponding to $\psi_1$
and $\psi_2$ are both isomorphisms, then the variety
${\mathcal A}_1$ is isomorphic to ${\mathcal A}_2$.

Therefore, if $X$ and $Y$ are two compact
connected Riemann
surfaces such that ${\rm Pic}^0(X)$ is isomorphic
to ${\rm Pic}^0(Y)$, then the moduli space of
rank one regular connections on $X$ is isomorphic to the
moduli space of rank one regular connections on $Y$. Since
there are non--isomorphic compact Riemann surfaces with
isomorphic Jacobians, we conclude that Theorem \ref{thm000}
fails for $r\, =\,1$.}
\end{remark}

\noindent
\textbf{Acknowledgements.}\, The first--named author is
grateful to N. Fakhruddin for useful discussions. We thank the
referee for pointing out a couple of issues to be clarified in
the first version of the paper. The second--named author is
partially supported through grant from MEC (Spain)
MTM2004-07090-C03-01.


\end{document}